\documentclass[twocolumn, amssymb, nobibnotes, aps, amsmath, amssymb, superscriptaddress,floatfix]{revtex4-2}
\usepackage{graphicx} 
\usepackage{xcolor}
\begin{document}

\title{Topologically Directed Simulations Reveal the Impact of Geometric Constraints on Knotted Proteins}
\author{Agnese Barbensi}
\affiliation{School of Mathematics and Physics, University of Queensland, Brisbane}
\author{Alexander R. Klotz}
\affiliation{Department of Physics and Astronomy, California State University, Long Beach}
\author{Dimos Gkountaroulis}
\affiliation{Department of Molecular and Human Genetics, Baylor College of Medicine, Houston}
\affiliation{The Center for Genome Architecture, Baylor College of Medicine, Houston}
\date{March 2025}

\begin{abstract}
    Simulations of knotting and unknotting in polymers or other filaments rely on random processes to facilitate topological changes. Here we introduce a method of \textit{topological steering} to determine the optimal pathway by which a filament may knot or unknot while subject to a given set of physics. The method involves measuring the knotoid spectrum of a space curve projected onto many surfaces and computing the mean unravelling number of those projections. Several perturbations of a curve can be generated stochastically, e.g. using the Langevin equation or crankshaft moves, and a gradient can be followed that maximizes or minimizes the topological complexity. We apply this method to a polymer model based on a growing self-avoiding tangent-sphere chain, which can be made to model proteins by imposing a constraint that the bending and twisting angles between successive spheres must maintain the distribution found in naturally occurring protein structures. We show that without these protein-like geometric constraints, topologically optimised polymers typically form alternating torus knots and composites thereof, similar to the stochastic knots predicted for long DNA. However, when the geometric constraints are imposed on the system, the frequency of twist knots increases, similar to the observed abundance of twist knots in protein structures.
\end{abstract}

\maketitle

\section{Introduction}
Polymers are frequently knotted, raising the question of whether these knots form and untie through \textit{preferred} mechanisms. For instance, in DNA molecules, complex knots tend to migrate towards the ends under elongational fields and usually untie in multiple stages, even when the unknotting number is 1~\cite{soh}. Conversely, the unfolding and folding of knotted proteins are often intricate, involving various intermediate states and multiple potential pathways (see \textit{e.g.}~\cite{mallam2006probing, lou2016knotted} or~\cite{tubiana2024topology} for a comprehensive review). The process by which a protein folds into its native state is in general highly complex, involving transitions through numerous locally minimal configurations within the folding landscape. Despite significant advances in protein structure prediction, particularly with the development of deep learning methods like AlphaFold~\cite{abramson2024accurate}, a comprehensive understanding of protein folding pathways remains out of reach. This is primarily because, although these methods provide accurate predictions, they offer little to no insight into the actual folding mechanisms. Knotted proteins, in addition to reaching a final folded configuration, also achieve a specific knotted topology, adding an extra layer of complexity to the folding process. However, knots also offer a perspective that simplifies the analysis of protein folding. By focusing on the topological changes of the folding intermediates, one can adopt a coarse-grained analysis of the pathway, rather than striving for a complete understanding of the transitions through local minimal states. In this sense, predicting a protein’s knotting pathway can serve as a simplified test case, providing insight into the more complex problem of protein folding pathway prediction. 

In this context, the search for \textit{preferred} (un)knotting mechanisms can be rephrased as follows: does the tying and untying of knots in biopolymers follow a \textit{topologically} optimal process? Or, said differently, does the \textit{complexity} of a biopolymer only change strictly monotonically along these pathways? In this project, we address the problem of knot folding by investigating this question. The main idea is to define a method of \textit{topological steering} based on topological measures of \textit{knot complexity}.

\begin{figure*}
    \centering
    \includegraphics[width=0.7\linewidth]{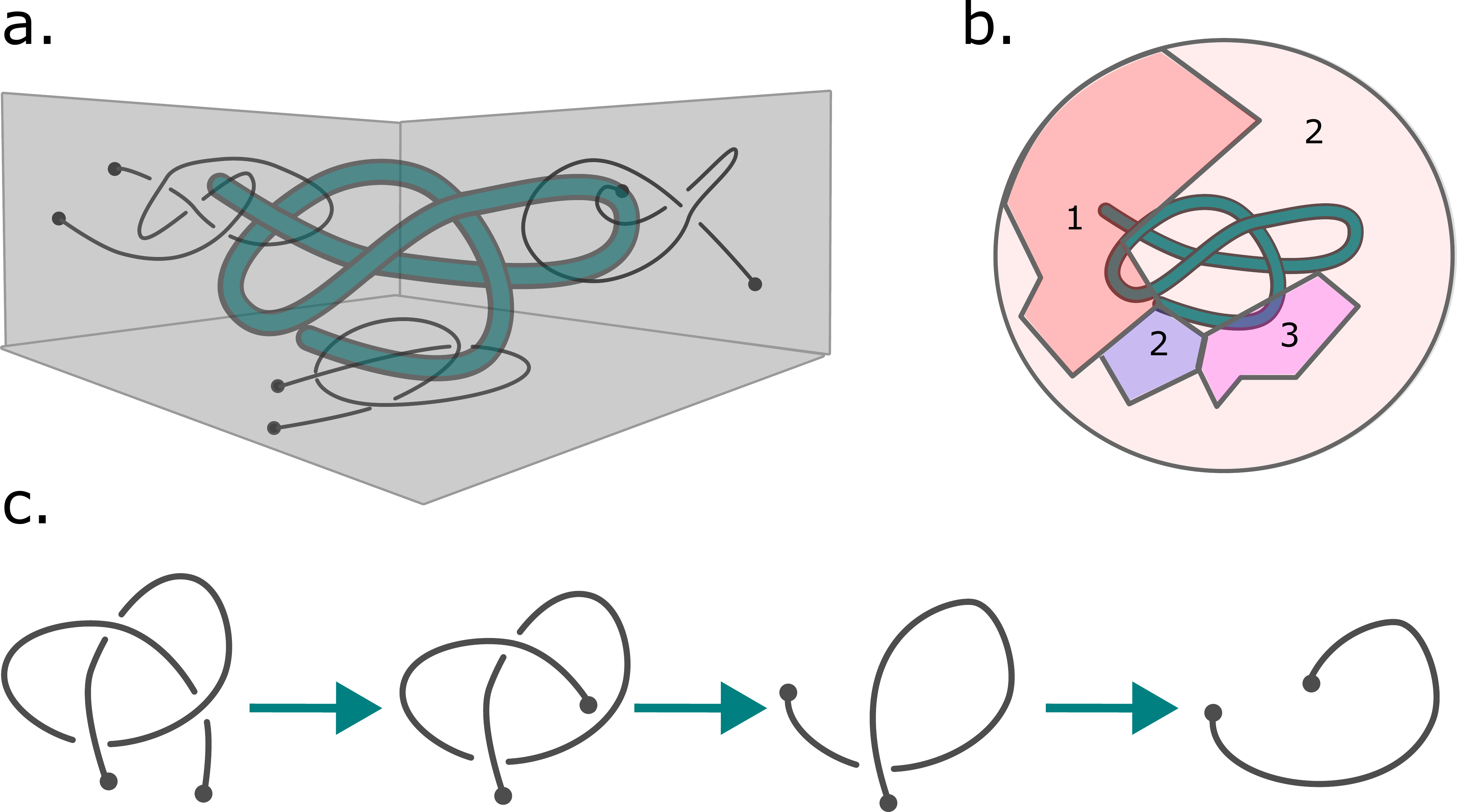}
    \caption{a. Projection of an open 3D curve onto three surfaces. Each projection corresponds to a knotoid type. b. A schematic for the $S^2$-distribution of knotoids of the curve in a. Each colour indicates a specific knotoid type, with the corresponding unravelling number added as a label. c. An example of a diagram unravelling into simpler knotoids through forbidden moves. The unravelling number of the left-most knotoid is equal to 2, as the penultimate diagram in the sequence differ by the trivial knotoid also by a single Reidemeister I move.} 
    \label{fig:1}
\end{figure*}

An essential step in this direction is being able to effectively describe the topology of knotted biopolymer, and how they change during the process. This is often done in the language of knot theory, which formally only applies to closed loops. Previous work has described optimal paths to unknotting or unlinking, but these typically focus on transitions between closed curves (e.g. due to magnetic reconnection or topoisomerase) \cite{liu2016knots, liu2020minimal}. Several schemes have been developed to assign a knot topology to an open space curve, to analyse experimental data or simulations of polymers. Among the most commonly used approaches are those that artificially \textit{close} the space curve to form a closed loop. These closure techniques, such as end-to-end, minimally interfering, and stochastic closure, characterise the space curve by the knot type it assumes upon closure~\cite{tubiana2011probing}. In the case of stochastic closure, the characterisation is provided by a distribution of knot types, rather than a single one. A different point of view consists in analysing the entanglement \textit{intrinsically}, \textit{i.e.}~without relying on knot closures. Intrinsic approaches include Vassiliev measures~\cite{panagiotou2021vassiliev, klotz2024} and knotoids~\cite{goundaroulis2017topological}, see also~\cite{tubiana2024topology} for a recent review. Similarly to stochastic closures, the knotoid approach describes the topology of an open curve through the distribution of its planar projections in the two-dimensional sphere $S^2$. These projections are intrinsically analysed as well-defined topological objects~\cite{goundaroulis2017topological}, see Figure~\ref{fig:1}(ab). The resulting knotoid distribution simultaneously captures the \textit{global} topology and the \textit{local} geometry of the entangled curve~\cite{dorier2018knoto, dabrowski2019knotprot, barbensi2021f, barbensi2021topological}.

The next key step is to use these knot-theoretical descriptions to define a meaningful notion of topological complexity. In this work, we develop a measure of the knot complexity of an open curve based on the average unravelling number~\cite{barbensi2021f} of knotoid projections of the curve, see Figure~\ref{fig:1}(c). This measure of complexity defines a \textit{topological functional} on the space of open curves; we then search for optimal (un)knotting pathways by developing an optimisation algorithm to maximise and minimise this functional, based on a gradient descent of many perturbations of the curve. To demonstrate the utility of this topological steering algorithm, we apply it to Langevin dynamics simulations of semiflexible polymers, showing how knotted polymers may optimally untie and how unknotted polymers may optimally entangle themselves.

The ability to direct simulated polymers towards greater topological complexity enables us to access an ensemble of knotted polymers, even when the knotting probability is low. Only about 1\% of known proteins contain knots~\cite{dabrowski2019knotprot, tubiana2024topology}. Biological proteins, selected by evolution, represent a tiny subset of all possible peptides up to a certain length. It remains unknown whether the proportion of knotted biological proteins is less than, greater than, or the same as the proportion of knots in all possible proteins. Without this knowledge, arguments about the biological advantages or disadvantages of knots in proteins are incomplete. Despite recent advancements in protein structure prediction, addressing this question computationally remains unfeasible. However, our topological steering protocol allows us to generate large samples of knotted protein-like polymers. In turn, this establishes a reference framework to explore both the mechanisms of knot formation and the specific types of knots that emerge.

In random knots, despite variations due to specific model choices, there is general agreement on the relative likelihood of different knot types. Simpler knots are more common than complex ones, and compound simple knots, rather than knots with higher crossing number, dominate as chain length increases \cite{virnau}. This pattern might suggest that a similar distribution of knot types will be observed in biopolymers. Proteins however are known to exhibit a strong preference for twist knots, appearing with a higher-than-random occurrence in AlphaFold-predicted structures~\cite{niemyska2022alphaknot, rubach2024alphaknot}, and with non-twist knots entirely absent in experimentally solved structures~\cite{dabrowski2019knotprot}. We propose that this phenomenon can be partially explained by the fact that proteins are subject to strong geometric constraints~\cite{prior2020obtaining}, which alter the landscape of achievable knots. To investigate this idea, we use the mathematical formalism of topologically optimal pathways and topological steering we developed. After establishing the validity of topological steering, we present a model of topologically-directed protein-like random walks, and we show that imposing the same geometric constraints found in nature for proteins results in a substantial increase in the formation of twist knots. This provides the first evidence that local polymer geometries are responsible for the protein knotting spectrum observed in nature.

\section{Main}
To explore the existence of topologically optimal (un)knotting pathways, we first introduce in Subsection~\ref{sec:AUN} the average unravelling number $\text{AUN}$, a measure of topological complexity for knotted arcs, defining a topological functional. We then explore how to approximate the corresponding gradient-flow induced dynamics via an optimisation method. To this end, we first consider the case of  molecular dynamics simulations of semiflexible Kremer-Grest polymers in Subsection~\ref{sec:semiflexible}. We then explore the space of configurations of the polymer by applying random Brownian force to each bead, and we move towards those configurations minimising (respectively maximising) $\text{AUN}$. We finally explore the effect of protein-like geometric constraints on the topological pathways and on the resulting knotted population using a growing self-avoiding walk model, described in Subsection~\ref{sec:GSW}

\subsection{Knotoids, the Average Unravelling Number and the Total Unravelling Number}~\label{sec:AUN} Knotoids are topological objects defined as equivalence classes of arc diagrams, subject to Reidemeister moves performed away from the endpoints, and to planar isotopies of $S^2$. Examples of knotoid diagrams are shown in Figure~\ref{fig:1}(ac). Knotoids have been used extensively to characterise the topology of knotted 3D arcs~\cite{goundaroulis2017topological,dorier2018knoto, dabrowski2019knotprot, barbensi2021f, barbensi2021topological}. Indeed, any generic planar projection of such a curve defines a knotoid diagram (see Figure~\ref{fig:1}(a)), whose knotoid type can be computed using topological invariants. By considering all possible planar projections of the curve and their corresponding knotoid types, one can obtain an $S^2$-distribution of knotoids which describes the arc's entanglement, as illustrated in Figure~\ref{fig:1}(b).

Any knotoid diagram can be transformed into a diagram for the trivial knotoid by a sequence of \textit{forbidden} moves, see Figure~\ref{fig:1}(c). Similarly to the case of knots, crossing changes and unknotting numbers, one can define the \textit{unravelling number}~\cite{barbensi2021f} as the minimal number of forbidden moves needed to transform a knotoid into the trivial one, among all the possible diagrams in the same class.

Consider a 3D curve defined by an embedding $\psi: [0,1] \longrightarrow \mathbb{R}^3$. We define the average unravelling number $\text{AUN}$ as the following integral: $$ \text{AUN} = \int_{S^2} u(k)\; \mathrm{d}S^2, $$ where the knotoids $k$ vary in the $S^2$-knotoid distribution of the curve $\psi([0,1])$. The AUN provides a measure of complexity for open 3D curves, and a natural candidate for a topological functional. 

\begin{figure}[h!]
    \centering
    \includegraphics[width=0.7\linewidth]{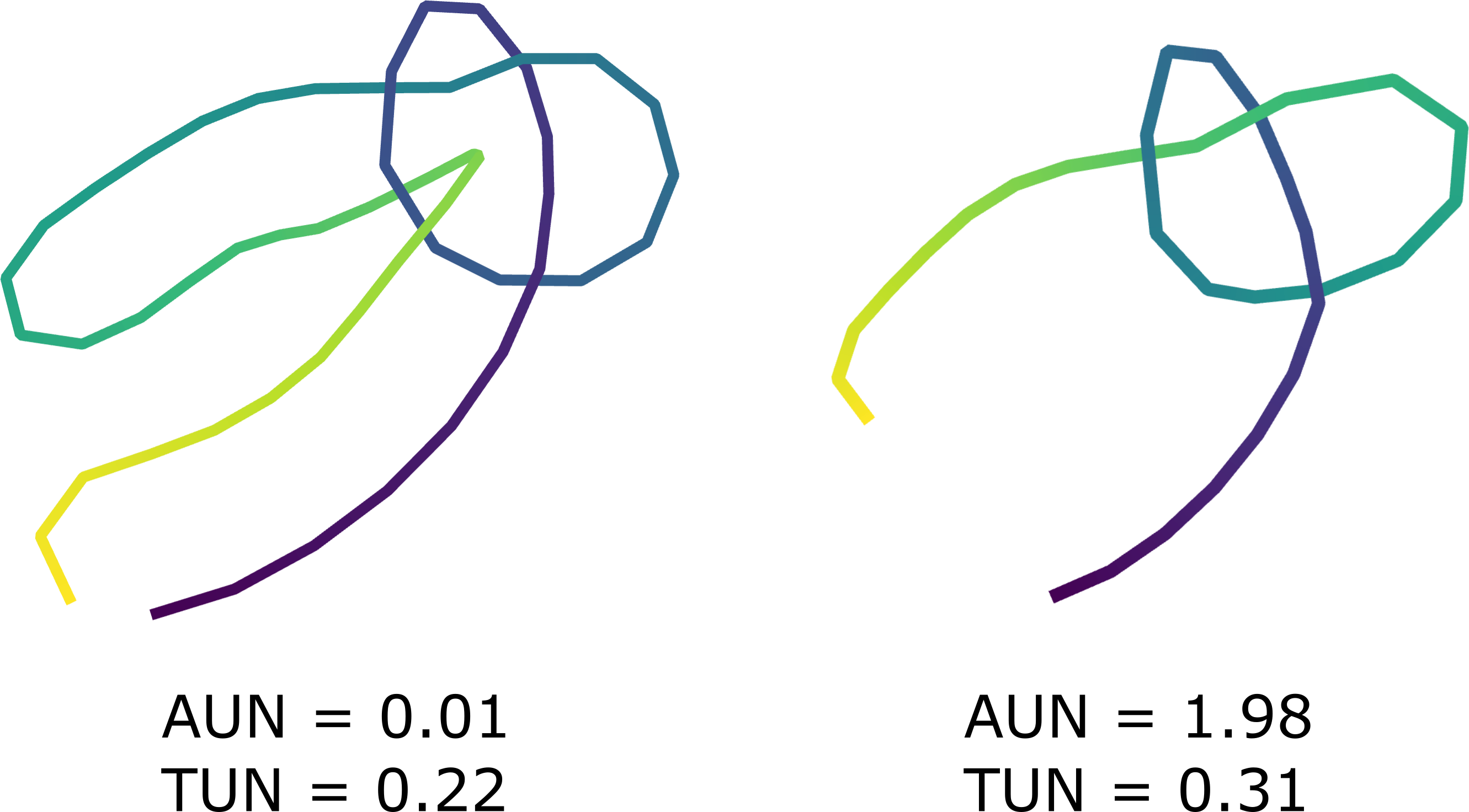}
    \caption{A 3D curve creating a slipknot (left), and an open-ended trefoil (right). For the slipknot, while the $\text{TUN}$ adequately captures the local entanglement of the curve, the $\text{AUN}$ is very close to $0$. }
    \label{fig:2}
\end{figure}

For curves creating \textit{slipknots}, however, the AUN tends to be close to zero (see Figure~\ref{fig:2}), which can result in local minima during gradient descent unknotting. To mitigate this issue, we define the Total Unravelling Number (TUN) as the integral of AUN over all subchains of the 3D curve: $$ \text{TUN} = \int_{\Delta} \text{AUN}(s)\; \mathrm{d}s. $$ Here, $s$ represents the pair $s = (a,b)$, which determines the subchain $\psi([a,b])$, where $a$ and $b$ vary within the triangle $\Delta$ such that $0 \leq a \leq b \leq 1$. This approach accounts for the contributions of all subchains, thereby avoiding the formation of local (slip)knotted minima and providing a more comprehensive measure of the curve's global entanglement.

Computationally, the knotoid distribution of a 3D curve and its corresponding AUN can be approximated by \textit{e.g.}~uniformly sampling points in $S^2$ to select finitely many projections. In the case of polygonal curves, we can then compute an approximated TUN by iteratively trimming points from both ends of the curve. 
In this work, we use the software Knoto-ID~\cite{dorier2018knoto} to approximate the knotoid distribution, while values for the unravelling number of knotoids are taken from~\cite{barbensi2021f}. For details on the computations, see the Software and Data availability statement. 
\subsection{Topologically Directed Simulations}\label{sec:semiflexible}

Ideally, given a knotted polymer, one would want to follow the AUN- or TUN-induced gradient-descent flow to produce unknotting pathways. However, differentiating these two functionals is currently out of reach. For this reason, we opt to explore the space of configurations of open knots in different physical models, and then apply optimisation procedures to move towards those conformations minimising topological complexity. As a proof of concept, we first we apply the topological optimisation procedure to molecular dynamics simulations of semiflexible Kremer-Grest polymers. The polymer is treated as a chain of beads connected to their neighbours by finitely extensible springs, with Lennard-Jones repulsion between non-bonded beads. The maximum spring extension is set to 1.5 times the Lennard-Jones diameter of the beads, making strand crossings extremely unlikely. Semiflexibility is imposed by a Kratky-Porod potential on the cosine of the angle between three successive beads, and is tuned such that the persistence length is 10 times the bead diameter. The configuration of the polymer is iterated forward in time by the Langevin equation, which imposes viscous drag and a random Brownian force on each bead. The equation of motion is solved in LAMMPS. This model has been used previously for knotted DNA-inspired polymers \cite{sleiman, klotz2024}.

\begin{figure*}
    \centering
    \includegraphics[width=0.98\linewidth]{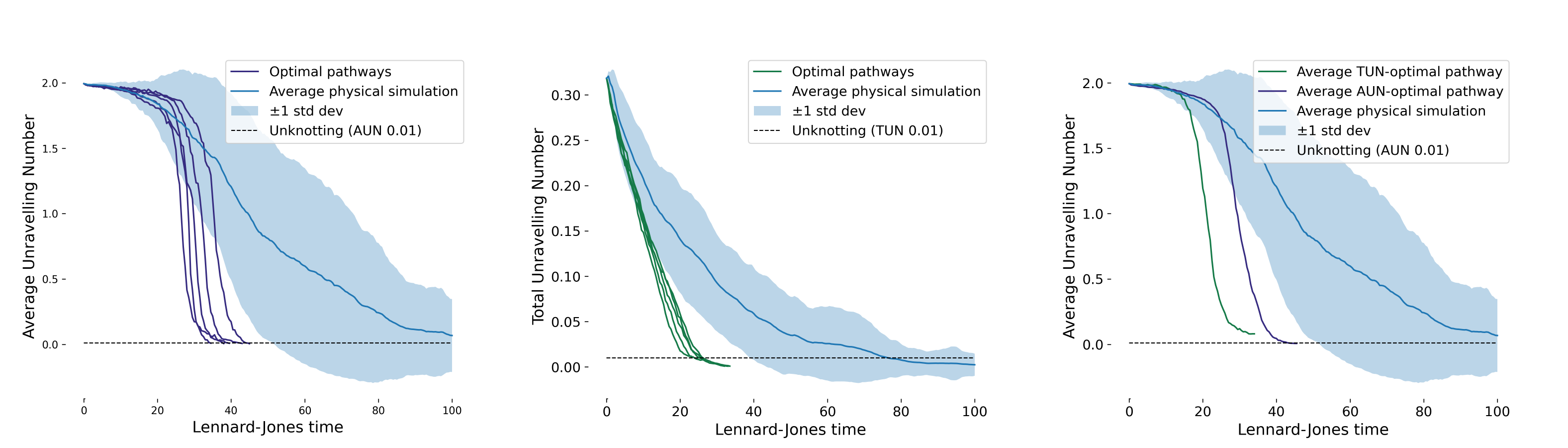}
    \caption{a. AUN of a semiflexible polymer initially in a deep $3_1$ knotted configuration as it unties. Shades of purple show the topologically optimal untying, while blue shows an ensemble average of several untying events driven by thermal fluctuations. b. Same as a., with TUN instead of AUN. Shades of green show the topologically optimal untying, while blue shows an ensemble average of several untying events driven by thermal fluctuations. c. Comparison between AUN and TUN-driven unknotting.} 
    \label{fig:unknotting}
\end{figure*}

We first consider simulations driven by the AUN. We interface our topological steering protocol with LAMMPS by taking an equilibrated initial configuration of length $38$, and simulating 40 independent runs with different random seeds, for an interval of 500 Langevin steps of 0.001 $\tau_{LJ}$ (where $\tau_{LJ}$ is the Lennard-Jones timescale defined the Supplemental Information). We then compute the $\text{AUN}$ for the final configuration of each parallel run, and we select the configuration with the minimum value as the new initial condition for the next iteration, with 40 new random seeds. The sequence of optimal configurations approximates the most efficient way for the topology of the polymer to evolve under the physics of the Kremer-Grest chain. We then replace AUN with TUN, and repeat the simulations. Videos of topologically optimised unknotting using both functionals are included as supplementary material, see the Software and Data availability statement. A sufficiently short open chain with a tight knot far from its ends will untie stochastically, as there are more unknotted than knotted configurations. The time to achieve unknotting has a broad distribution, as studied \textit{e.g.}~by Lai \cite{lai}, Carraglio et al. \cite{caraglio}, and Klotz and Estabrooks \cite{klotz2024}. Figure \ref{fig:unknotting}(ab) shows AUN and TUN of a topologically-directed unknotting polymer as a function of the number of iterations. In comparison, the population average of many undirected untying polymers are shown. We observe that in both cases, the directed simulation reaches a baseline unknotted state faster than the stochastic population, and --after perhaps an initial steadiness-- its topological complexity falls linearly rather than exponentially, indicating consistent progress towards the unknot. A comparison between AUN and TUN-driven trajectories shows that the latter unknot faster, see Figure~\ref{fig:unknotting}. As shown by the unknotting videos and by the evolution of the topological complexities in Figure~\ref{fig:unknotting}(ab), however, the optimal unknotting trajectories appear qualitatively similar independently on the choice of AUN vs TUN. In both cases, the videos show the arms of the chain progressively moving into the knot as the knot expands, rather than the arms fluctuating randomly as typically occurs with polymer knots.

\begin{figure}[h!]
    \centering
    \includegraphics[width=0.98\linewidth]{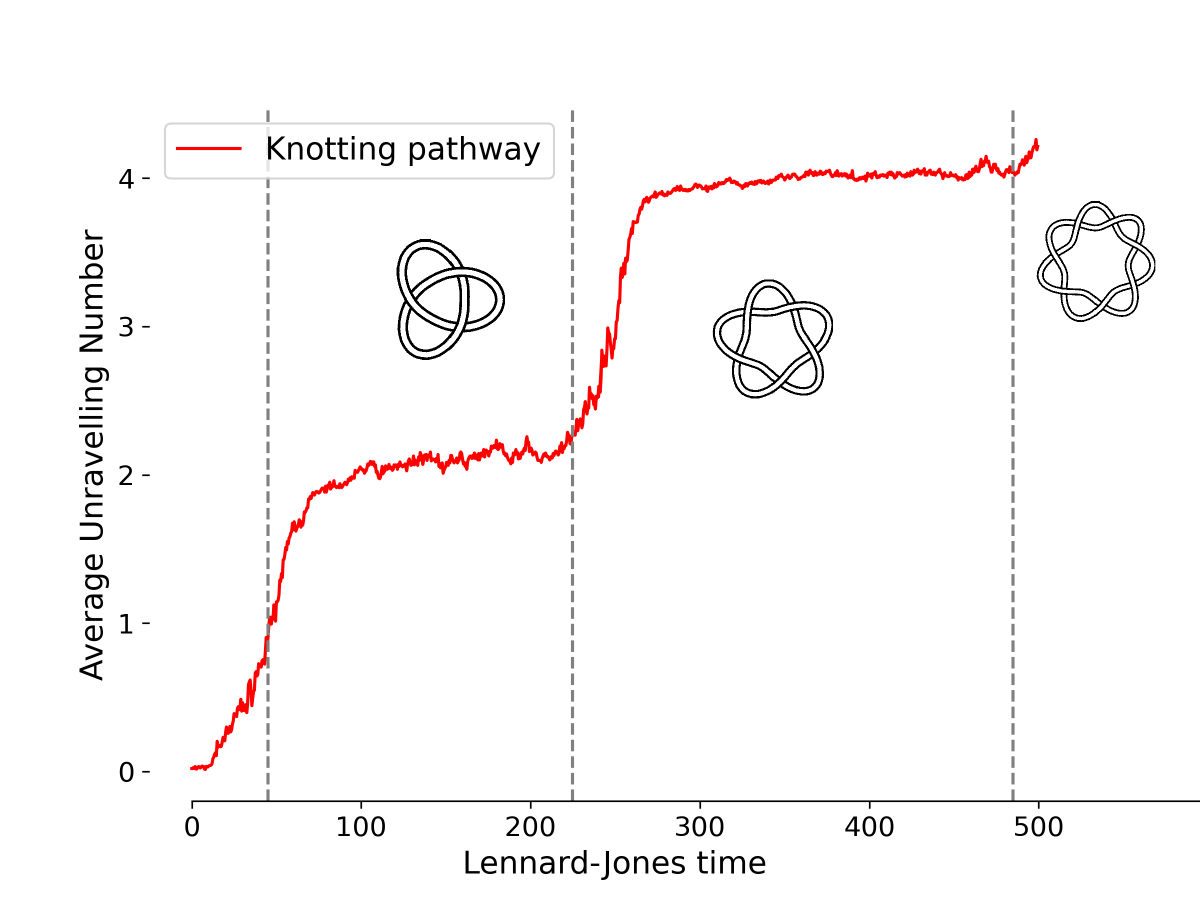}
    \caption{Topologically optimal knotting of an initially unknotted open semiflexible chain, forming alternating torus knots of increasing complexity.}
    \label{fig:knotting}
\end{figure}

We next apply the method to investigate optimal knotting pathways, focusing on AUN as computing this quantity is much more feasible compared to TUN. The polymers in our simulations are typically too short (length 52) to spontaneously knot, and a stochastic sampling of uncorrelated states will yield curves with non-trivial AUN approximately 1\% of the time, while the probability of actual knots is below 0.05\% (see the supplementary information and the Software and Data availability statement for more details). When the simulations are steered to maximise (rather than minimise) the knotoid complexity, we observe that the chain appears initially stochastic until it is able to form a trefoil knot, then increases its complexity by having the ends of the chain wrap around the existing core of the knot. The complexity of the open polymer steers the simulations that ultimately determine the knot type of the underlying curve, which is determined via stocastic closures using the software Knoto-ID~\cite{dorier2018knoto}. This leads to a sequence of alternating torus knots of increasing complexity ($3_1$, $5_1$, $7_1$...) as the simulation proceeds, see Figure~\ref{fig:knotting}. A possible explanation for this behaviour is that the simulated polymers generally tend to spontaneously form twisted loops as their initial configurations with non-zero AUN. Once these loops are formed, repeatedly piercing them becomes the quickest way to increase complexity, given how the functional AUN is defined. This raises a natural question: could the introduction of geometric constraints alter the spectrum of initial non-trivial configurations that spontaneously form, and consequently, affect the resulting knot population?

\subsection{Growing Walk Model}\label{sec:GSW} 
\begin{figure}[h!]
    \centering
    \includegraphics[width=0.5\linewidth]{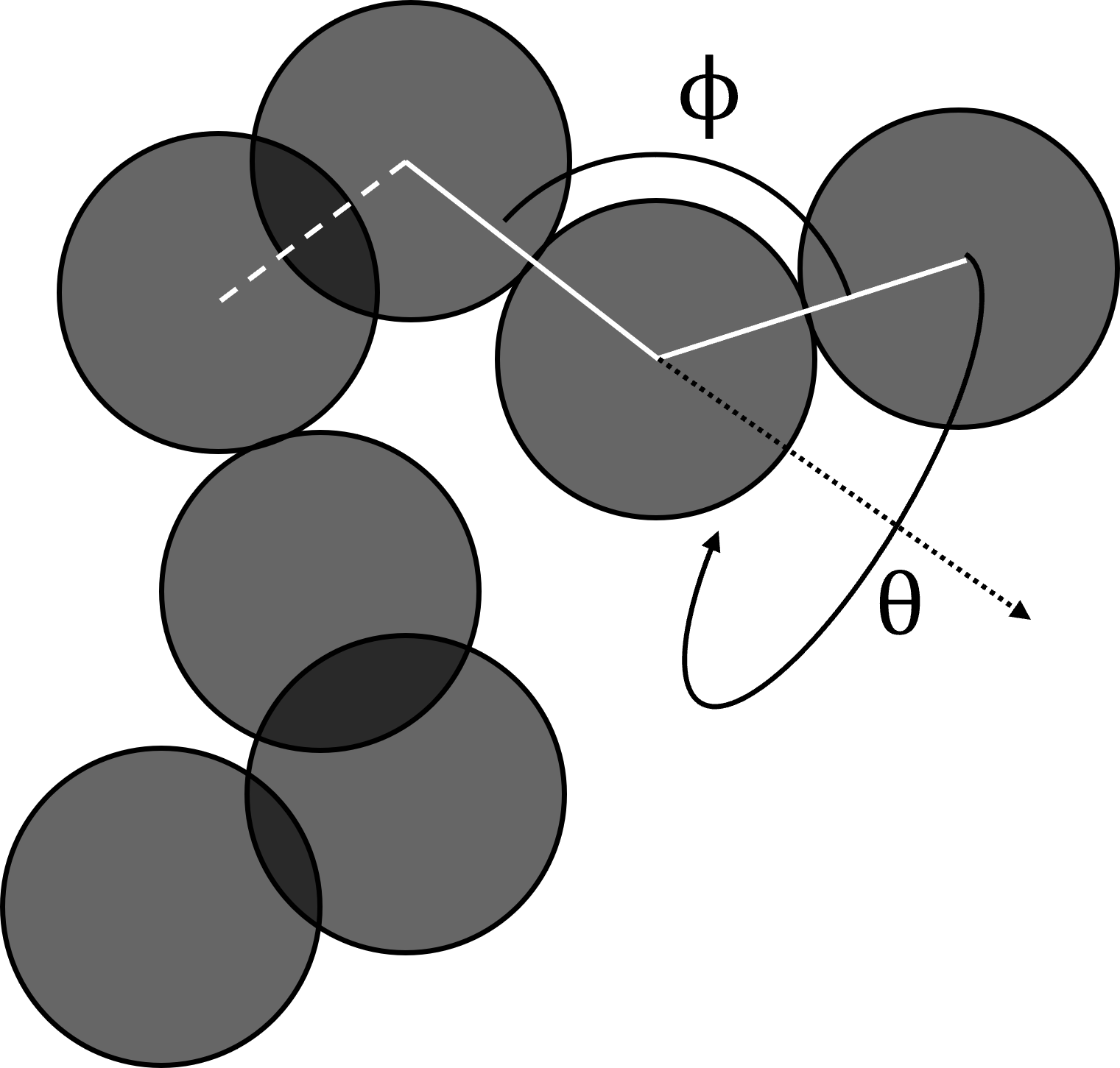}
    \caption{Schematic of the growing tangent-sphere chain, in which successive spheres are added according to a bending angle $\theta$ and a dihedral angle $\phi$.}
    \label{fig:angles}
\end{figure}

To facilitate the inclusion of local geometric constraints, we implement a growing self-avoiding walk (GSAW) model. In this model, we treat monomers as beads of unit radius and grow a chain by placing the $n$-\emph{th} bead tangent to the $(n-1)$-\emph{th} without overlapping any bead in the chain. The local geometry of a polymer can be described in terms of the bending angle $\theta$ and an azimuthal dihedral angle $\phi$, see Figure~\ref{fig:angles}. These angles in turn determine the local curvature $\kappa$ and torsion $\tau$.

Inspired by the method of~\cite{prior2020obtaining}, in our GSAW model the new position is chosen according to the angles $\phi$ and $\theta$. This allows us to impose different geometric constraints. Indeed, when both angles are chosen according to a uniform distribution, a self-avoiding freely-jointed chain is grown. If instead then angle $\phi$ is chosen from a Gaussian distribution about $\pi$, the walk gains a persistence length inversely proportional to the variance of the distribution. In practice, to prevent self-intersections, angles are chosen from a given distribution and the new coordinate is checked for overlaps. If any overlaps are present, then new angles are chosen until a bead is successfully placed, or a preset number of attempts is reached. If no bead is found then the walk may be terminated, similar to trapping self-avoiding walks on lattices \cite{klotz2020}, or a bead may be placed in an overlapping position making it ``weakly self-avoiding'' similar to the Domb-Joyce model \cite{dombjoyce}. We note that  however, even for a flexible walk, trapping is highly unlikely below length 1000, which is longer than the typical lengths we study under topological steering. For more information on the computational details, see the Software and Data availability statement.

\begin{figure}[h!]
    \centering
    \includegraphics[width=1\linewidth]{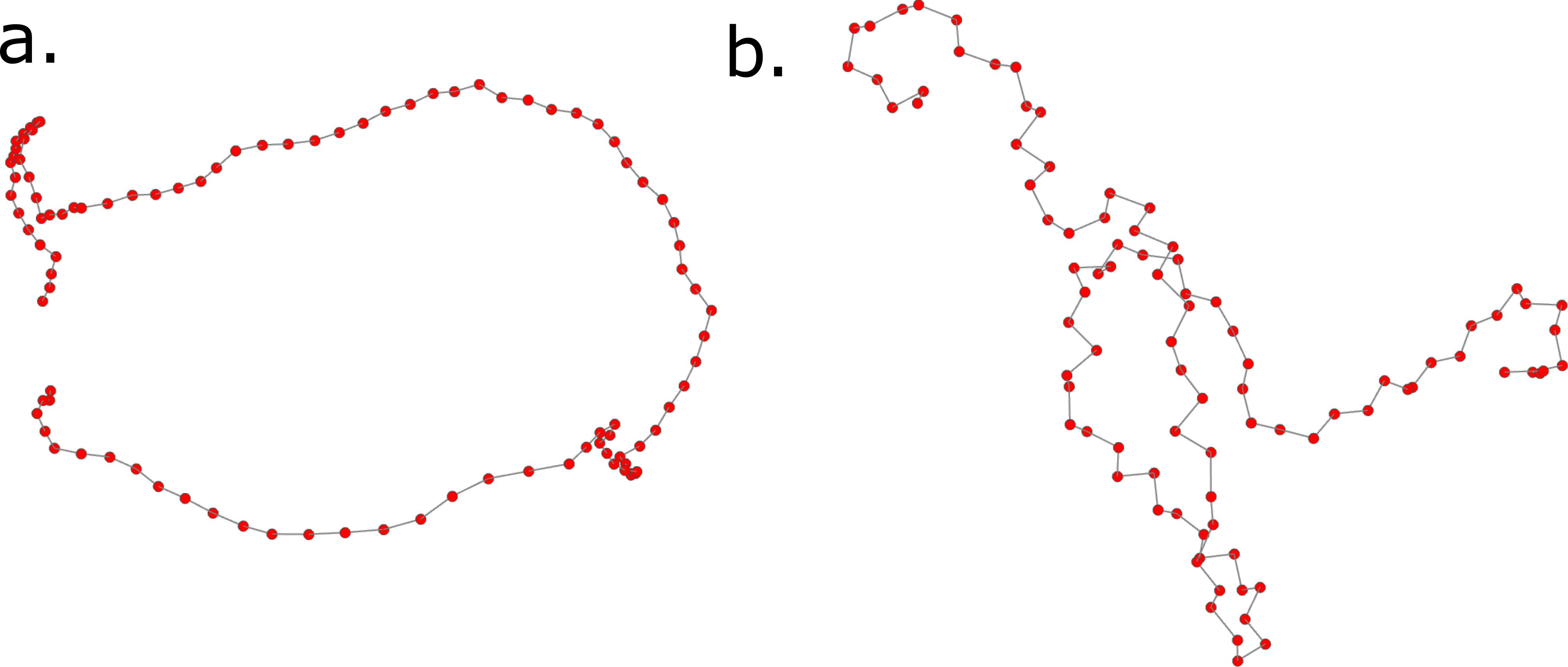}
    \caption{a. Example of an unbiased walk of length 100. b. Example of a protein-like walk of length 100.}
    \label{fig:unbiased}
\end{figure}

Protein structures exhibit a very specific local geometry, with highly skewed curvature $\kappa$ and torsion $\tau$ distributions~\cite{prior2020obtaining}. To replicate this behaviour, we modify the model by treating each monomer as a C$\alpha$ atom in a polypeptide, and by drawing $\phi$ and $\theta$ from the distribution found in native protein structures.

\begin{figure}[h!]
    \centering
    \includegraphics[width=0.9\linewidth]{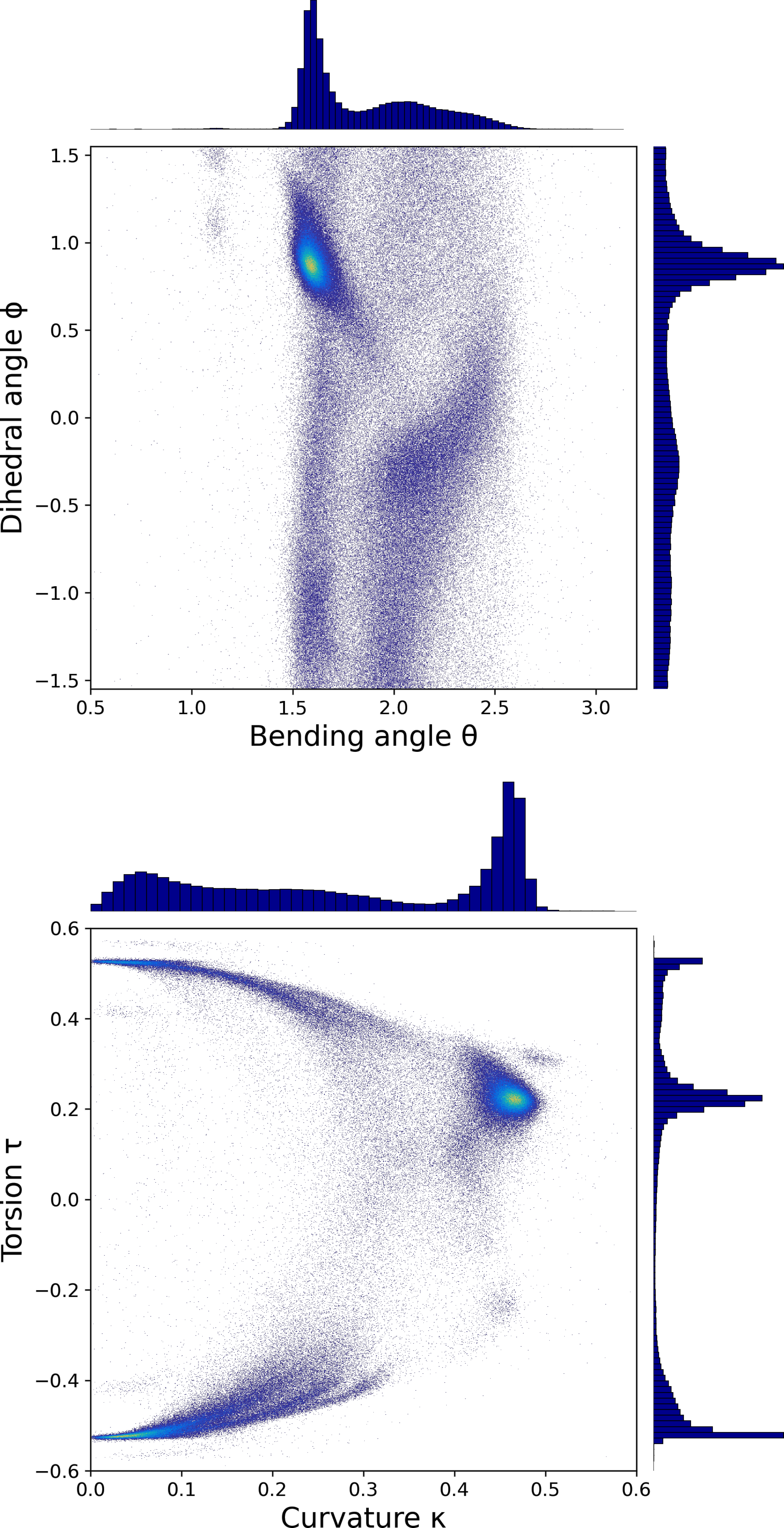}
    \caption{Distributions bending angle $\theta$ and a dihedral angle $\phi$ (top), and curvature $\kappa$ and torsion $\tau$ found in proteins.} 
    \label{fig:protangs}
\end{figure}

\begin{figure*}
    \centering
    \includegraphics[width=0.7\linewidth]{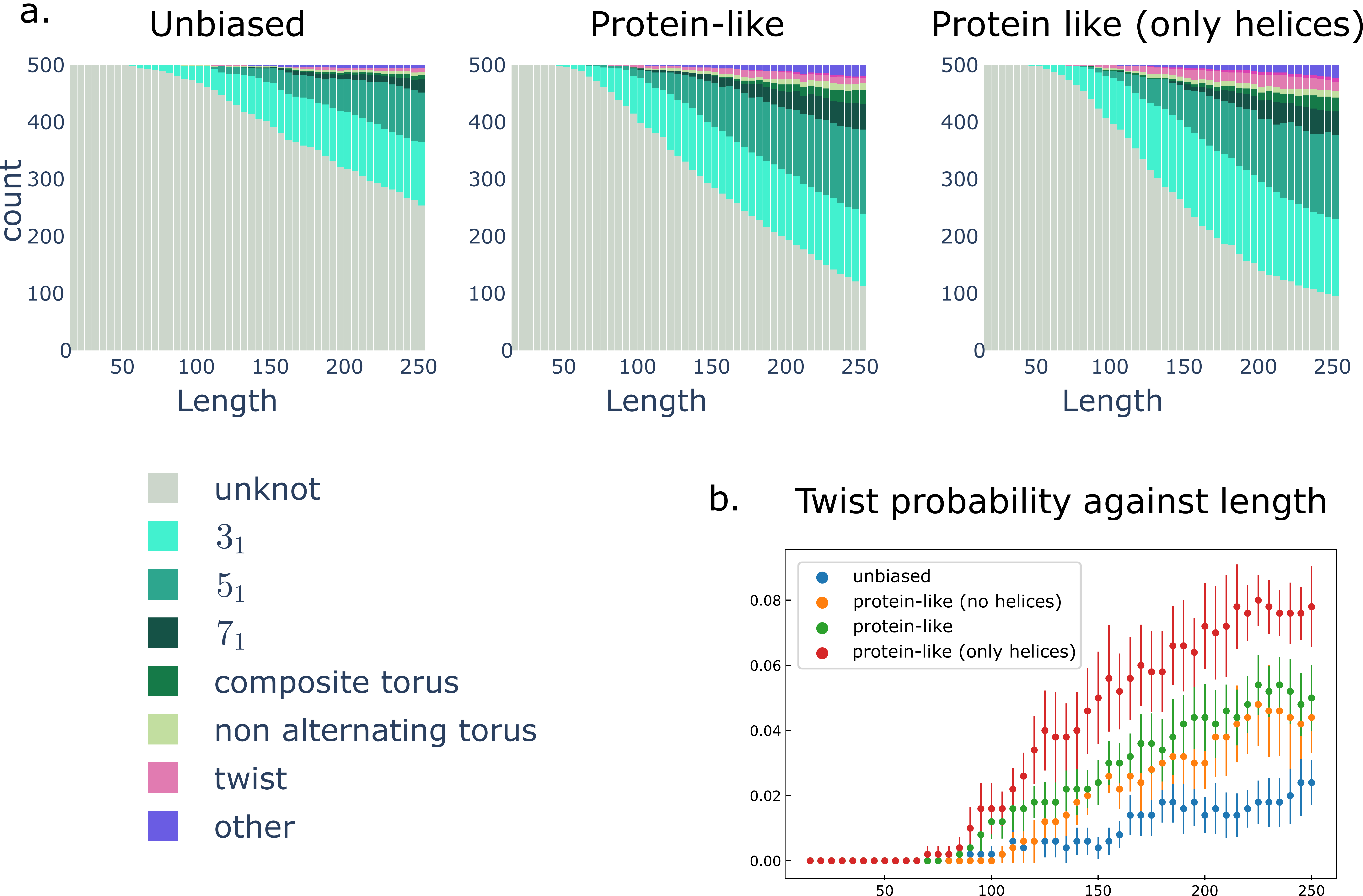}
    \caption{a. Kymographs of knot formation during topologically-steered walk growth, showing the fraction of different knot types in a population of 500 walks as a function of length, for the different geometric models considered. For unbiased walks, near 50\% of the knots are trefoils. As the protein-like geometry is implemented, we observe an increase in diversity of the knot population and a higher probability of occurrence of twist knots. b. Proportion of twist knots as a function of length for topologically directed growing walks, for four different models. Error bars are measured over random subsamples of the data. Note that knot types here are determined via stocastic closures, using the software Knoto-ID~\cite{dorier2018knoto}.}
    \label{fig:agres}
\end{figure*}

To do so, we measure the bending angles between every three successive beads and the dihedral angles between every four successive beads in 1000 crystallographically-determined non-redudant protein structures (see \textit{e.g.}~\cite{prior2020obtaining}), yielding 179,617 correlated $\phi$-$\theta$ pairs (see Figure~\ref{fig:protangs}), and curvature and torsion distributions replicating those found in proteins. The correlated angle distribution has two peaks, corresponding to helical and non-helical sections of the proteins. We grow \textit{protein-like} random walks by drawing an angle pair at random and using it to affix the new bead to the chain. This produces random walks that have the $\phi$-$\theta$ distribution as well as the curvature-torsion distribution found in native proteins, see Figure~\ref{fig:unbiased}(b). We can additionally constrain the choice of angles to draw exclusively from the helical or non-helical section of the distribution. We call these walks \textit{protein like (no helices)} and \textit{protein-like (only helices)}, respectively. As a point of comparison, we generate semiflexible random walks with no dihedral constraint and a persistence length of 10 monomers, comparable to DNA in aqueous solution, as this is the model typically used to study knot formation in random biopolymers \cite{virnau}. We call these walks \textit{unbiased}, see Figure~\ref{fig:unbiased}(a).

\section{Results and Discussion}
Similar to the implementation of Section~\ref{sec:semiflexible}, we can steer the growth of these polymers by adding a bias towards configurations maximising AUN. In practice, in each of the unbiased, protein-like, protein-like (no helices) and protein-like (only helices) models, we grow the chain randomly by adding 10 beads in 20 different ways, we compute the AUN of all the resulting chains, and we pick the growth maximising the complexity. Figure~\ref{fig:agres}(a) shows the kymographs of ensembles of these topologically steered growing walks. For each length reached by the population of walks, the colour indicates the fraction of independent runs with a given knot type. The first observation is that for each of the models, the knotting probability exceeds 50\% by length 250, which is quite high compared to the undirected probability of 2-4\% (see the supplementary information). We also observe that the protein-like models have a notably higher knotting probability, with the maximum values reached by drawing angles exclusively from the helical section of the angle distribution. 
In all the models, trefoil knots are by far the most common, and the phenomenon of creating increasingly complex alternating torus knots observed in Section~\ref{sec:semiflexible} is essentially replicated. However, as the protein-like constraints are added, we see that the fraction of alternating torus knots is substantially higher in the unbiased sample than in the protein-like models. In contrast, the protein-like populations displays a higher proportion of twist knots, as well as more complex knots beyond the twist/torus families. We note that these differences, although small, are statistically significant, as demonstrated by Figure~\ref{fig:agres}(b), presenting the proportion of twist knots in each ensemble as a function of the length. The only-helices protein-like walks present the highest proportion, reaching 8\% at length 250, while the unbiased walks have the lowest, reaching about 2\%. The same bias towards twist knots can be observed in undirected growing models, although, the very low knot probability makes the differences substantially harder to detect and to statistically validate, see the Supplementary Information for more details.

These simulations were motivated by the question of why knotted proteins have a higher proportion of twist knots than would be expected in stochastic polymer knots. Our hypothesis is that this discrepancy arises from the geometric constraints imposed on proteins by their bending and twisting angles. A recent analysis of over 600,000 knotted protein structures, as determined by AlphaFold \cite{niemyska2022alphaknot, rubach2024alphaknot}, revealed that approximately 9.3\% were non-trefoil twist knots. For knots with experimentally determined crystal structures, this ratio increases to 13.4\%~\cite{dabrowski2019knotprot}. Notably, there are over 12 times as many $5_2$ (twist) knotted proteins as $5_1$ (torus) knotted proteins. For protein knots with 5, 6, and 7 crossings, the twist variant is the most common. Note that the typical number of amino acids in a knotted protein is around 250, which matches the terminal length used in our simulations. In contrast, simulations of knotting in DNA show that as the length increases, the chains tend to form higher-order composite trefoils.

These results provide the first compelling evidence that geometric angle constraints play a crucial role in shaping the protein knotting spectrum observed for proteins in nature.

\section{Conclusion}
In this study, we have investigated the utility of topologically steered simulations in studying polymer (un)knotting. We introduce the average unravelling number (AUN) as a measure of topological complexity for knotted arcs. By defining a topological functional and approximating the corresponding gradient-flow induced dynamics, we have demonstrated the potential of topological steering in guiding polymers towards configurations of varying complexity. Our simulations, revealed that topological steering can effectively direct polymers towards both knotted and unknotted states. This approach not only enhances our understanding of the mechanisms underlying knot formation and untying, but also provides a powerful tool for generating large samples of knotted protein-like polymers. Such samples are very effective for studying the formation and types of knots in biological proteins, which remain a topic of significant interest. While our findings strongly suggest that geometric angle constraints play a crucial role in the protein knotting spectrum observed in nature, the individual knot probabilities still fall short of accurately replicating those found in proteins. Modifications to the protein model more reminiscent of those used to study folding, including interactions between amino acids their ellipsoidal geometry, will be a fruitful path towards understanding protein knots when coupled to topological steering. 

\section*{}

\subsection*{Software and Data availability statement}
The software used in this paper, and the data generated, are freely available on Zenodo~\url{https://zenodo.org/records/15220263}. Our computations rely on Knoto-ID~\cite{dorier2018knoto}. You can either compile it from source or download a pre-built executable from its GitHub repository.

\subsection*{Acknowledgements}
The authors would like to thank Christopher Prior, Arron Bale, and Raffaele Potestio for helpful comments. A substantial part of this work was carried out during a visit of AK to AB, kindly supported by the Raybould Visitor Fellowship, School of Mathematics and Statistics, the University of Queensland. AK is additionally supported by the National Science Foundation, grant number 2336744.

\end{document}